\documentclass{article}
\usepackage{amsmath}
\usepackage{amscd}
\usepackage{amsthm}
\usepackage{amssymb} \usepackage{latexsym}
\usepackage{eufrak}
\usepackage{euscript}
\usepackage{epsfig}
\usepackage{graphics}
\usepackage{caption}
\usepackage{subcaption}
\usepackage{array} 
\usepackage{enumerate}

\usepackage{tikz}

\usepackage{color}

\newcommand{\red}[1]{\textcolor{red}{#1}}

\usepackage{hyperref}

\newcommand{\bel}[1]{\begin{equation}\label{#1}}

\newcommand{\be}{\begin{equation}}

\newcommand{\ba}{\begin{eqnarray}}
\newcommand{\ea}{\end{eqnarray}}
\newcommand{\rf}[1]{(\ref{#1})}

\newcommand{\qe}{\end{equation}}

\newcommand{\R}{{\mathbb R}}

\newtheorem{thesis}{Thesis}
\newcommand{\btl}[1]{\begin{thesis}\label{#1}}
\newcommand{\et}{\end{thesis}}

\theoremstyle{theorem}
\newtheorem{theo}{Theorem}[section]
\newtheorem{prin}{Principle}[section]

\theoremstyle{corollary}

\theoremstyle{lemma}
\newtheorem{lemma}{Lemma}[section]
\theoremstyle{definition}
\newtheorem{defi}{Definition}[section]
\theoremstyle{proof}

\theoremstyle{remark}

\title{Geometry of Data}
\author{Parvaneh Joharinad and J\"urgen Jost}
\begin{document}
\maketitle

\begin{abstract}
Topological data analysis asks when balls in a metric space $(X,d)$
intersect. Geometric data analysis asks how much balls have to be enlarged to
intersect.\\
We connect this principle to the traditional core geometric concept of
curvature. This enables us, on one hand, to reconceptualize curvature and link
it to the geometric notion of hyperconvexity. On the other hand, we can then
also understand methods of topological data analysis from a geometric perspective.
\end{abstract}

\section{Introduction}
Many data sets come with a basic geometric structure, distances between data
points. It is therefore natural to use geometric methods to analyze such
data. The deepest geometric concepts, however, were developed in the 19th
century for smooth
manifolds, more precisely Riemannian manifolds. And the most fundamental
concept there is curvature. In the 20th century, notions of curvature were
successfully generalized to more general classes of spaces. Still, those
spaces, like geodesic length spaces,  are typically not discrete, in contrast
to data sets. Thus, we have found it desirable to rethink fundamental
geometric concepts from a more abstract perspective that also naturally
includes discrete spaces. Of course, there are ideas and approaches that we
can build upon, most importantly those pioneered by Gromov \cite{Gromov80,Gromov99}. From such a
perspective, the distinction between discrete and connected spaces is partly
one of scale. From a large scale perspective, spaces from those two classes
may look alike. 

Such a large scale perspective is still quantitative, hence geometric, and is
therefore different from a qualitative topological approach. Nevertheless,
as we shall see, there are important links between the two. In particular, we
can look at the successful topological data analysis method of persistent
homology from a geometric perspective. 

Topological data analysis asks when balls in a metric space $(X,d)$
intersect. This is a qualitative concept, but the data analysis method of
persistent homology makes this quantitative through the dependence on the
radii of the balls. Geometric data analysis, as we conceive it in this
contribution,  asks how much balls have to be enlarged to intersect. And as we
shall see, this is captured by a suitable concept of curvature. And curvature,
from a general perspective as adopted here, quantifies convexity. Therefore,
convexity and its strengthening as hyperconvexity will be our basic concepts.

\section{Preliminaries from metric geometry}
Let $(X,d)$ be a  metric space. $x,y,\dots $ will be points in $X$, and they
thus have a distance $d(x,y)$. A continuous  path $c:[0,1]\longrightarrow X$  with  $x=c(0), y=c(1)$ has length
$$\mathit{l}(c):=\sup\sum_{i=1}^{i=n}d(c(t_i),c(t_{i-1})).$$ 
The supremum  here  is taken over  all partitions of $[0,1]$, with $t_0=0,
  t_{n}=1$. $(X,d)$ is called a \emph{length space} if  for all $x,y$,
$$d(x,y)=\inf\{\mathit{l}(c):\: \text{c is a path between}\; x \; \text{and}\;
y\}.$$
A length space $(X,d)$ is  called {\em geodesic} if this infimum is always
realized, that is, any $x,y\in X$ can be  connected by
a  {\em shortest  path} $c:[0,1]\longrightarrow X$,
i.e.
$$d(x,y)=\mathit{l}(c).$$
Thus, the distance between $x$ and $y$  is realized by some curve, a shortest
geodesic.\\
Every complete locally compact length space is a geodesic space. However, there is another way to determine whether a complete metric space is a geodesic (resp. length) space by checking the existence of mid-points (resp. approximate midpoints).
\begin{defi}
  $m\in X$ is  a {\em midpoint} between $x,y$  if
  $$d(x,m)=d(m,y)=\frac{1}{2}d(x,y).$$
\end{defi}
We may also say that a pair of points $x,y\in X$ has  {\em approximate midpoints} if for every $\epsilon>0$ there exists $m_\epsilon\in X$ with
\[
\max\{d(m_\epsilon,x),d(m_\epsilon,y)\}\leq\frac{1}{2}d(x,y)+\epsilon
\]
We observe 
\begin{lemma}
Every pair of points in a geodesic space (resp. length space) has at least one midpoint (resp. approximate midpoints).\\
The inverse is true provided that the metric space is complete.  \qed
\end{lemma}
In the sequel,
$$ B(x,r):=\{ y\in X: d(x,y)\le r\}$$
will always be the {\em closed} ball centered at $x$ with radius $r\ge 0$. 
\begin{defi}\label{totally convex}
  $(X,d)$ is     {\em totally convex} if 
  for any   $x_1,x_2\in X, r_1, r_2 >0$ with
  $${r_1+r_2\geq d(x_1,x_2)},$$
  we have 
  $$B(x_1,r_1)\cap B(x_2,r_2)\neq \emptyset .$$ 
\end{defi}
Any radii $r_i$ will be $>0$ in the sequel.\\
Again, an easy lemma
\begin{lemma}
  {Geodesic spaces} are totally convex.\qed
\end{lemma}
Length spaces are not necessarily totally convex, as they need not be
complete. An example is  $\R^2\setminus{0}$ with the length structure induced
by the Euclidean distance.\\

Let us formulate Definition \ref{totally convex} as a
\begin{prin}
  {\bf Two balls that can intersect do intersect.}
\end{prin}

We shall now introduce a fundamental quantity. For
$$r_1+r_2\geq d(x_1,x_2)$$
we put
\begin{eqnarray}
  \label{10a}
  \rho((x_1,x_2),(r_1,r_2))&:=&\inf_{x\in X}\max_{i=1,2}\frac{d(x_i,x)}{r_i}\\
  \label{10b}
\rho(x_1,x_2)&:=&\sup_{r_1,r_2} \rho((x_1,x_2),(r_1,r_2))
\end{eqnarray}
If $\rho(x_1,x_2)=1$ for each pair of points $x_1,x_2\in X$, then the
existence of approximate midpoints is guaranteed, and  $X$ is a length space
provided that it is a complete metric space. If, moreover, the infimum is
attained for each pair by some $x_0\in X$, then $X$ is a geodesic space
provided that it is complete.\\
Another obvious
\begin{lemma}
When $X$ is complete the supremum in (\ref{10b}) is realized by $r_1=r_2=\frac{1}{2}d(x_1,x_2)$, that is
\begin{equation}\label{11}
\rho(x_1,x_2)=\inf_{x\in X}\max_{i=1,2}\frac{2d(x_i,x)}{d(x_1,x_2)}.
\end{equation}
Moreover, $\rho(x_1,x_2)=1$ is achieved for some $x$ when
$$d(x_1,x)+d(x_2,x)=d(x_1,x_2),$$
that is, when $x$ is a midpoint of $x_1,x_2$. 
\end{lemma} \qed

Thus, we want to find points {between} two points $x_1$ and $x_2$, and quantify to
  what extent that can fail.\\
 Therefore, in the realm of complete metric spaces, the more (\ref{10b})
 deviates from $1$ the less is the chance to approximate distances by lengths
 of connecting paths.\\

  A key idea now is to extend this to three points.

   \section{Tripod spaces}
\begin{defi}
   A geodesic length space $(X,d)$ is a {\em tripod space} if for any three points $x_1,x_2,x_3 \in X$, there exists a {\em median,} that is, a point $m\in X$ with
\begin{equation*}
d(x_i,m)+d(x_j,m)=d(x_i,x_j),\: \text{for} \: 1\leq i<j\leq 3.
\end{equation*}
\end{defi}
We note that for a median, we have
$$d(x_1,x_2)+d(x_2,x_3)+d(x_3,x_1)=2(d(x_1,m)+d(x_2,m)+d(x_3,m)).$$
\begin{center}
\begin{tikzpicture}[scale =1]
 \node[above] (a) at (1,6) {$x_1$};
  \node[above] (b) at (5,6) {$x_2$};
  \node[right] (c) at (3,1) {$x_3$};
  \node[right] (d) at (3,4) {\red{$m$}};
  \fill (1,6) circle (3pt);
  \fill (5,6) circle (3pt);
  \fill (3,1) circle (3pt);
  \fill[red] (3,4) circle (3pt);
  \draw (1,6) -- (3,4);
  \draw (5,6) -- (3,4);
  \draw (3,1) -- (3,4);
\end{tikzpicture} 
\end{center}

Most metric spaces are {\em not} tripod spaces. For instance,
  Riemannian manifolds of dimension $>1$ do {\em not} satisfy tripod
  property. Nevertheless, there are {\bf examples} that will be important for
  us:
  \begin{itemize}
\item Metric trees
\item $L^\infty$-spaces
\item  and more generally, {\em hyperconvex spaces} (to be defined shortly)
\end{itemize}
If such a median exists it will be a minimizer for the sum of the distances to
the corresponding triple $x_1,x_2,x_3$. Such a point is  called a Fermat point.\\

 {\bf Our strategy} will then be to quantify the  deviation from the  tripod  property. \\

 We get the existence of tripods if the following more general condition is
 satisfied. 
For any  $x_1,x_2,x_3 \in X$ which  do not lie on a geodesic, and  $r_i+r_j\geq d(x_i,x_j)$, $1\leq i<j\leq3$, 
$$\bigcap\limits_{i=1}^{3}B(x_i,r_i)\neq \emptyset.$$
This leads to
\begin{prin}
 {\bf Three balls that can intersect do intersect.} 
\end{prin}
To explore this principle, and the deviation from it, we shall now introduce a 3-point analogue of \rf{10a}, \rf{10b}
For $x_1,x_2,x_3 \in X$ and $r_i+r_j\geq d(x_i,x_j)$, 
\begin{eqnarray}
  \label{20a}
  \rho((x_1,x_2,x_3),(r_1,r_2,r_3))&:=&\inf_{x\in
                                        X}\max_{i=1,2,3}\frac{d(x_i,x)}{r_i}\\
  \label{20b}
\rho(x_1,x_2,x_3)&:=&\sup_{r_i+r_j\ge  d(x_i,x_j), i\neq j} \rho((x_1,x_2,x_3),(r_1,r_2,r_3)).
\end{eqnarray}
This is uniquely solved by the {\em Gromov products} 
\begin{align}\label{Gromov}
\nonumber r_1&=\frac{1}{2}(d(x_1,x_2)+d(x_1,x_3)-d(x_2,x_3)),&\\
\nonumber r_2&=\frac{1}{2}(d(x_1,x_2)+d(x_2,x_3)-d(x_1,x_3)),&\\
r_3&=\frac{1}{2}(d(x_1,x_3)+d(x_2,x_3)-d(x_1,x_2)).&
\end{align}
\textbf{Remark:} It is obvious that $\rho((x_1,x_2,x_3),(r_1,r_2,r_3))\geq 1$. Moreover, this quantity is bounded from above by $2$ if $X$ is  complete.\\
If (with $r_1,r_2,r_3$ defined by (\ref{Gromov})) $\rho(x_1,x_2,x_3)=1$ and the infimum is attained by some $m$, then we have a tripod construction or equivalently a Fermat point. This implies that there exists an intermediate point through which each pair $x_i,x_j$ can be connected. 
\begin{defi}
 An $m$ attaining the  infimum in (\ref{20a}) is called  a {\em weighted circumcenter}. 
\end{defi}
A weighted circumcenter solves an optimization problem in $\R^3$ with respect
to the $\mathit{l}_{\infty}$ norm. The larger the value of $\rho(x_1,x_2,x_3)$ is, the less optimal the weighted circumcenter as the interconnecting point will be.\\
We observe here
\begin{lemma}
  Weighted circumcenters exist and are unique  for  triangles in $CAT(0)$ spaces 
(Alexandrov's generalization of Riemannian manifolds of sectional curvature  $\le 0$).
\end{lemma}

\section{Hyperconvexity}
\begin{figure}
\centering
\begin{subfigure}{.5\textwidth}
  \centering
  \includegraphics[width=.8\linewidth]{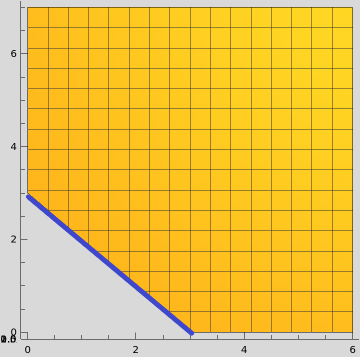}
  \caption{2-point}
  \label{2-point}
\end{subfigure}%
\begin{subfigure}{.5\textwidth}
  \centering
  \includegraphics[width=.8\linewidth]{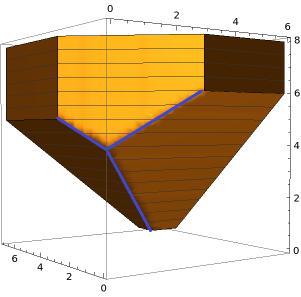}
  \caption{3-point}
  \label{3-point}
\end{subfigure}
\caption{a) The yellow area is the set of all possible radius functions on $2$ points and the line segment colored in blue refers to the minimal ones.\\
b) The three dimensional polyhedron is the set of all possible radius functions on $3$ points and the tripod consisting of three line segments colored in blue refers to the minimal ones.}
\label{fig:test}
\end{figure}

We shall now extend the above principle to arbitrary numbers of points. 
\begin{defi}\label{hyperconvex}
 $(X,d)$ is  {\em hyperconvex} if 
 for any family $\{x_i\}_{i\in I}\subset X$ and  {$r_i+r_j\geq d(x_i,x_j)$} for $i,j\in I$, 
\[
\bigcap\limits_{i\in I}B(x_i,r_i) \neq \emptyset
\]
\end{defi}
In a totally convex metric space, $r_i+r_j\geq d(x_i,x_j)$ can be replaced by $
B(x_i,r_i)\cap B(x_j,r_j)\neq\emptyset$ for all $i,j\in I$. Thus, when
balls intersect pairwise, they also have a common intersection.\\
This leads to our final 
  \begin{prin}
{\bf Balls that can intersect do intersect.}
 \end{prin}
 We observe
 \begin{lemma}
  Hyperconvex spaces are tripod spaces. \qed
 \end{lemma}
 We list some important properties of hyperconvex spaces
 \begin{theo}
  \begin{enumerate}
\item[a)] Hyperconvex spaces  are complete and  contractible to each of their
  points \cite{Aronszajn56}.
\item[b)]  $X$ is hyperconvex iff every $1-$Lipschitz map from a subspace of
  any metric space $Y$ to $X$ can be extended to a $1-$Lipschitz map over $Y$
   \cite{Aronszajn56}.
\item[c)] Every metric space is isometrically embedded in a hyperconvex space, called its hyperconvex hull. The hyperconvex hull of a compact space is compact and that of a finite space is a simplicial complex. \cite{Isbell64, Dress84}
\end{enumerate} 
\end{theo}
We now describe the isometric embedding in the part (c) and the
  construction of the hyperconvex hull, in order to understand the specific
  choice of radii in (\ref{10a}) and (\ref{20a}).  By the Kuratowski embedding, every metric space $(X,d)$ is isometrically embedded in the space of bounded functions on $X$ equipped with the supremum norm, i.e. $l_{\infty}(X)$, via the map $x\mapsto d(x,.)$ which we denote by $x\mapsto d_x$ for simplicity.\\
$l_{\infty}(X)$ contains the subspace $E(X)$ consisting of all functions $f$ that are minimal subject to the  relation
\bel{hyp1}
f(x)+f(y)\geq d(x,y),\; \forall x,y\in X.
\qe
It has been shown in \cite{Isbell64, Dress84, Lang13} that $E(X)$ is a
hyperconvex space containing the image of  $X$ under the Kuratowski embedding
isometrically, and  $E(X)$  is minimal in the sense that it is isometrically embeddable in any other such  hyperconvex space.\\
The radii in (\ref{10a}) and (\ref{20a}) are functions
on a $2$-point space and a $3$-point spaces respectively, satisfying (\ref{hyp1}). \\
If $X$ is a finite metric space with $|X|=n$, the space of all functions
satisfying (\ref{hyp1}) is a polyhedron in the finite vector space $\R^n$
obtained by the intersection of the closed half spaces $f_i+f_j\geq
d(x_i,x_j)$ for $1\leq i<j\leq n$. Therefore,  the  interior of every face $S$
of this polyhedron is the intersection of some hyperplanes
$f_i+f_j=d(x_i,x_j)$. We can then define a graph $G(S)$ with vertex set $X$,
corresponding to the symmetric relation defined by that face. More precisely,
$x_i$ is connected to $x_j$ with an edge in $G(S)$ if for $f\in S^\circ$ we
have $f_i+f_j=d(x_i,x_j)$. Now, $E(X)$ is the union of compact faces of this
polyhedron and moreover the graph corresponding to each such face is a
spanning graph, that is, every vertex is connected to at least one other
vertex in this graph. This construction was first introduced in
\cite{Dress84}, where a combinatorial dimension for finite metric spaces was
defined as the maximal dimension of a face in its  hyperconvex hull. The
hyperconvex hull of finite metric spaces was  studied further in
\cite{Bandelt92, Dress02} and from a different perspective in \cite{DeLoera95, Sturmfels04} to obtain the \textit{metric fan} of a finite set. In \cite{Sturmfels04} a software tool was presented to visualize these hyperconvex hulls. The problem of finding faces of $E(X)$, when $X$ is finite, as a linear programming problem was also studied in \cite{Hirai06,Develin06} \\  
In the special case $X=\{x_1,x_2\}$ with distance $d_{12}=d(x_1,x_2)$, the corresponding polyhedron is the half plane $f_1+f_2\geq d_{12}$ cut by the coordinate planes $f_i=0,\; i=1,2$,  which has only one compact face, the line segment $f_1+f_r=d_{12}$ connecting $(d_{12},0)$ to $(0,d_{12})$, i.e., $[d_{x_2},d_{x_1}]$. Every point in this polyhedron can be reached through a ray passing this line segment and the midpoint of this segment, that is $\frac{1}{2}(d_{12},d_{12})$ is the corresponding radius function in (\ref{11}). \ The space of all such radius functions is illustrated in Figure \ref{2-point}.\\
Similarly, one can see that for $X=\{x_1,x_2,x_3\}$, using the same notation
$d_{ij},\; 1\leq i,j\leq 3$ for pairwise distances, the corresponding
polyhedron is the intersection of the half-spaces 
\be
f_i+f_j\geq d_{ij},\: 1\leq i< j\leq 3
\qe
and  the coordinate half spaces $f_i\geq 0$ for $i=1,2,3$. Moreover, the hyperconvex hull, colored in blue in Figure \ref{3-point}, is the union of three segments each of which connect a distance function $d_{x_i}$ to the function $r=(r_1,r_2,r_3)$ defined in (\ref{Gromov}). 
For the analysis of discrete metric spaces, some variants of the notion of
hyperconvexity are well suited, c.f \cite{Khamsi00, Espinola01,  Grunbaum59, Grunbaum60}. 
 \begin{defi}
  $(X,d)$ is   $\delta$-hyperbolic ($\delta \ge 0$) if  for any family $\{B(x_i,r_i)\}_{i\in I}$   with $r_i+r_j\geq d(x_i,x_j)$, 
\bel{del1}
\bigcap_{i\in I}B(x_i,\delta+r_i)\neq\emptyset.
\qe 
\end{defi}
 \begin{defi}
  $(X,d)$ is    $\lambda$-hyperconvex ($\lambda\geq1$) if for every family $\{B(x_{i}, r_{i})\}_{i\in I}$ with  $r_{i}+ r_{j}\geq d(x_{i}, x_{j})$, 
\bel{del2}\bigcap\limits_{i\in I}B(x_{i},\lambda r_{i})\neq\emptyset.
\qe
\end{defi}
Of course, $0$-hyperbolicity and $1$-hyperconvexity are simply hyperconvexity. 
For large radii,  $\delta$   insignificant, and the concept of
$\delta$-hyperbolicity is therefore  good for asymptotic
considerations. In contrast, $\lambda$-hperconvexity  is invariant under
scaling the metric $d$, and it can therefore capture scaling invariant
properties of a metric space.\\
The preceding concepts allow for a quantification of the deviation from
hyperconvexity. The following results are known.
\begin{theo}
 Hilbert spaces are $\sqrt{2}$-hyperconvex. Reflexive  and dual Banach spaces
 are $2$-hyperconvex. Therefore, for  a measure space $(X,\mu)$,  $L^p(X,\mu), 1<p<\infty$, are $2$-hyperconvex, and if $X$ is finite,  $L^1(X,\mu)$  is also $2$-hyperconvex. $L^\infty(X,\mu)$ is hyperconvex. \cite{Khamsi00, Espinola01}  
\end{theo}

 \section{Relation with Topological Data Analysis (TDA)}
 \begin{defi}
For  a  family $(x_i)_{i\in I}$ in a metric space $(X, d)$ and $r>0$, we define the  \emph{\v{C}ech complex} $\check{C}_r((x_i), X)$ containing a $q$-simplex whenever
$$ \bigcap_{i=1,\dots ,q+1} B(x_i,r)\neq \emptyset .$$
{Here $(x_i)_{i\in I}$ is called the landmark set and $X$ is the witness
  set. When the  witness set coincides with the landmarks, we thus  define a non-empty intersection inside the sample set  $(x_i)_{i\in I}$ as the criterion for  a simplex.} We also define the \emph{Vietoris-Rips complex} $VR_r((x_i), X)$ containing a $q$-simplex whenever
 $$ B(x_i,r)\cap B(x_j,r) \neq \emptyset \ \text{ for all }i,j\in I.$$
\end{defi}

The two structures are not as different as they might appear, as the
difference between the criteria  for spanning a
simplex is whether the 
vertex set is contained in a ball of radius or of diameter $r$.\\
The principle of the important topological data analysis scheme of \emph{persistent homology} then is to record how the homology of these complexes varies as a function of $r$. \cite{Edelsbrunner00, Zomorodian05, Edelsbrunner08, Carlsson09}\\
Of course, every simplex of the \v{C}ech complex is also a simplex of the
Vietoris-Rips complex, but not necessarily conversely { unless for each
  simplex at least one of the balls of diameter $r$ containing the vertex set
  of that simplex has a center in the witness set.\\}
Deviation from hyperconvexity lets the  Vietoris-Rips complex contain more
simplices  than
the  \v{C}ech complex, or conversely
\begin{lemma}
  In a hyperconvex space, all  simplices that are filled in the Vietoris-Rips
  complex are also filled in the \v{C}ech complex. In particular, there is no
  contribution to local homology from unfilled simplices. \qed
\end{lemma}
{For instance, we can take a sample $(x_i)_{i\in I}$  from a geodesic metric
  space $(X,d)$ and compare $VR_r((x_i), X)$ with $\check{C}_r((x_i),
  E(X))$. For  the latter complex,  we take the hyperconvex hull of $X$, i.e.
  $E(X)$, as the witness set. It is clear that $\check{C}_r((x_i),
  E(X))\subset VR_r((x_i), X)$, as $X$ is a geodesic space and hence totally
  convex. Conversely, every simplex in $VR_r((x_i), X)$ is defined according
  to the criterion that balls of radius $r$ around its vertices intersect
  pairwise, which by hyperconvexity of $E(X)$ implies the existence of a
  common point between them in $E(X)$. In other words, the Vietoris-Rips
  complex of a metric family $(x_i)_{i\in I}, d$ coincides with its \v{C}ech
  complex but with different witness sets. This natural principle has been
  used in \cite{Lim21} to study the metric thickening of $S^1$ in its
  hyperconvex hull. A thorough study of the \v{C}ech and the  Vietoris-Rips
  filtration of $S^1$ can be found  in \cite{Adamaszek16,Adamaszek17}.\\
If $X$ is a closed Riemannian manifold, for small-enough radius $r$ depending
on the injectivity radius and a curvature bound, $VR_r(X)$ is homotopy
equivalent to $X$ by a well known theorem of Hausmann \cite{Hausmann95}. On
the other hand according to the nerve lemma, whenever $X$ is a paracompact
space and the family of open balls around sample points $(x_i)_{i\in I}$
with radius $r>0$ define a cover such that the non-empty intersections of any
finite number of them is contractible, the \v{C}ech complex
$\check{C}_{<r}((x_i), X)$ is homotopy equivalent to the original space $X$,
c.f. \cite{Hatcher01}. Although Hausmann's theorem is restricted to the case
where the original space, from which the sample is taken, is a Riemannian
manifold, both construction at some point reveal the topology of the
space. However, the Vietoris-Rips filtration ignores the geometry of the space
beyond the pairwise relations. The extent to which  higher order relations are
overlooked by considering Vietoris-Rips complexes can be quantified by
computing the deviation from hyperconvexity of different orders. This measures
how much one must expand balls to obtain a simplex in the \v{C}ech complex of
$(x_i)$ with witness set $X$ after that simplex is observed in the  \v{C}ech
complex of $(x_i)$ with witness set $E(X)$. The upper bound  $2$ for this
scale  is usually stated in the TDA literature, but this bound is not sharp.\\ 
For instance,  let us consider equilateral triangles  of perimeter $3a$ in the
Euclidean plane, in a  circle and in a metric tree.  That is,
$(x_1,x_2,x_3)$, $(x'_1,x'_2,x'_3)$ and $(\bar{x}_1,\bar{x}_2,\bar{x}_3)$ are
comparison triangles in the Euclidean plane, a circle and a hyperconvex space,
respectively. As noted in  (\ref{Gromov}), $r=\frac{a}{2}$ is the radius at
which each of these triples forms a simplex in the corresponding Vietoris-Rips
complex. However, we only need the upper bound of $2$ in the case of
$(x'_1,x'_2,x'_3)$, where the point are sampled from a circle which has the
highest deviation from hyperconvexity, for expanding the balls to obtain the
simplex in the \v{C}ech complex, c.f. \cite{Joharinad20}.\\   }
One can also more generally let the radii of the balls be different. That is,  for a vertex set $(x_i)_{i\in I}$ and a corresponding non-negative radius function $r$, we define the \v{C}ech complex containing a $q$-simplex $x_1,...,x_{q+1}$ whenever
$$ \bigcap_{i=1,\dots ,q+1} B(x_i,r(x_i))\neq \emptyset .$$
The Vietoris-Rips complex is defined in a similar way. And one can then look
at the resulting constructions for all such radius functions simultaneously  \cite{Joharinad20}.

\section{Curvature}
 We can use the preceding  concepts to compare spaces with each other, or with
 reference spaces, like Euclidean space. In geometry, such a comparison  is
 quantified by the concept of {\em curvature.}  From our abstract perspective, curvature relates intersection patterns of balls to
 convexity properties of distance functions. \\

 As pointed out by Klingenberg \cite{Klingenberg82}, the beginning of the
 theory of spaces of negative curvature can be dated to the work of von
 Mangoldt \cite{vonMangoldt81} in 1881 who showed that on a complete simply
 connected surface of negative curvature, geodesics starting at the same point
 diverge and can never meet again. This implies that the exponential map is a
 diffeomorphism. Apparently unaware of von Mangoldt's work,  Hadamard
 \cite{Hadamard1898} in 1898 proved further results about geodesics on
 surfaces of negative curvature. E.Cartan \cite{Cartan25} later considered
 negatively curved Riemannian manifolds of any dimension. For our purposes,
 non-positive, as opposed to negative, curvature is the appropriate concept,
 as we are interested in comparison theorems. 

  Let us first recall a by now classical concept of non-positive curvature, introduced by Alexandrov \cite{Alexandrov57}. 
 \begin{defi}
   The geodesic space $(X,d)$ is a {$CAT(0)$}{\em -space} if for all geodesics $c_1, c_2:[0,1]\longrightarrow X$ with $c_1(0)=c_2(0)$
 \begin{equation}\label{CAT(0)}
 d(c_1(t),c_2(s))\leq \|{\bar{c}_1(t)-\bar{c}_2(s)}\|,\, \forall\ t,s\in[0,1]
 \end{equation}
where $\bar{c}_1,\bar{c}_2:[0,1]\longrightarrow \mathbb{R}^2$ are the  sides
of the Euclidean comparison triangle in $\mathbb{R}^2$  with the same side
lengths as the triangle  $\bigtriangleup(c_1(0),c_1(1),c_2(1))$.
 \end{defi}
According to this definition, triangles in $CAT(0)$-spaces are not thicker than Euclidean triangles
with the same side lengths, c.f. \cite{Jost97, Burago01,Bridson13,  Alexander19}\\

There is another important concept of non-positive curvature, introduced by
Busemann \cite{Busemann48}.
\begin{defi}
  A geodesic space $(X,d)$ is a {\em Busemann convex space} if for every two geodesics $c_1, c_2:[0,1]\longrightarrow X$ with $c_1(0)=c_2(0)$, the distance function $t\mapsto d(c_1(t),c_2(t))$ is convex. 
\end{defi}
Geodesics in Busemann space diverge at least as fast as in Euclidean space.\\
Every $CAT(0)$ space is Busemann convex but not conversely. For complete Riemannian
manifolds, however, the two definitions agree and are equivalent to
non-positive sectional curvature in the sense of Riemann. \\
Several  generalizations of  these definitions to metric spaces that are not
necessarily  geodesic have been proposed, for instance  \cite{Bacak15, Berg18,
Alexander19}. 
We now present our definition from \cite{Joharinad19}.
\begin{defi}\label{npc}
  The metric space $(X,d)$ has \emph{non-positive curvature} if for each triple $(x_1,x_2,x_3)$ in $X$ with the comparison triangle  $\bigtriangleup(\bar{x}_1,\bar{x}_2,\bar{x}_3)$ in $\R^2$, one has
\begin{equation*}
\rho(x_1,x_2,x_3)\leq \rho(\bar{x}_1,\bar{x}_2,\bar{x}_3),
\end{equation*}
where $\rho(\bar{x}_1,\bar{x}_2,\bar{x}_3)$ is similarly defined by 
\be\nonumber
\rho(\bar{x}_1,\bar{x}_2,\bar{x}_3):=\min_{x\in \R^2}\max_{i=1,2,3}\dfrac{\|x-\bar{x}_i\|}{r_i}.
\qe
\end{defi}
According to this definition, the circumcenter of a triangle in a
non-positively curved space is at least as close to the vertices as in the Euclidean case. In other words, there is chance of finding a better intermediate point for each triple of points in such a space than in Euclidean plane.\\
For any triple of closed balls $\{B(x_{i}, r_{i});\:i=1,2,3\}$ with pairwise
intersection,  $\bigcap_{i=1,2,3}B(x_{i},\rho r_{i})$ is non-empty whenever
$B(\bar{x}_{i},\rho r_{i})$, $i=1,2,3$, have a common point. Thus, balls do
not need to be enlarged more than in Euclidean case to get triple
intersection. Thus, we can again formulate a 
\begin{prin}
 {\bf Balls intersect at least as easily as in Euclidean space.} 
\end{prin}
{\bf Examples:}
\begin{itemize}
\item 
Tripod spaces have non-positive curvature in the sense of
  Def. \ref{npc}, because there, $\rho=1$, which is the smallest possible value. 
\item Complete $CAT(0)$ spaces have non-positive curvature in the sense of
  Def. \ref{npc}. The converse not true; in fact, our spaces need not be geodesic, nor have unique geodesics.
\item {Approximate version applies to discrete spaces.} This is obviously
  important for questions of data analysis, and this in fact constitutes one
  of the motivations for Def. \ref{npc}.
\end{itemize}
We also have
\begin{theo}\label{npcgivessec}
A complete Riemannian manifold $(N,g)$ has non-positive curvature iff it has non-positive sectional curvature, c.f. \cite{Joharinad19}.
\end{theo}
Obviously, with the same concepts and constructions, one can also define upper
curvature bounds other than 0, by comparison with suitably scaled 2-spheres or
hyperbolic planes.

\section{Conclusions}
The \v{C}ech construction assigns  to a cover $\mathcal{U}=(U_i)_{i \in I}$ of
$X$ a simplicial complex $\Sigma(\mathcal{U})$ with vertex set $I$ and a
simplex $\sigma_J$ whenever $\bigcap_{j\in J} U_j\neq \emptyset$ for $J\subset I$.
 When all intersections are contractible, the homology of
$\Sigma(\mathcal{U})$ equals that of $X$ (under some rather general
topological conditions on $X$). When $(X,d)$ is  metric space, we can use
covers by (open or closed) distance balls. Now, when $(X,d)$ is a {\it hyperconvex} metric space, and if we use a cover $\mathcal{U}$ by distance balls, then whenever 
\bel{cc1}
\bigcap_{j\in J\backslash \{j_0\}} U_j\neq \emptyset \text{ for every }j_0\in J,
\qe
then also
\bel{cc2}
\bigcap_{j\in J} U_j\neq \emptyset,
\qe
 i.e., whenever $\Sigma(\mathcal{U})$ contains all the boundary facets of some
 simplex, it also contains that simplex itself. It even satisfies the stronger condition that whenever $\Sigma(\mathcal{U})$ contains all the boundary faces of dimension $1$ of some
 simplex, it also contains that simplex itself. This means that  $\Sigma(\mathcal{U})$ is a
 flag complex. Thus, there are no holes of the type of unfilled simplices, and no corresponding
 contributions to homology groups.\\

 As hyperconvex spaces are contractible, then whenever non-trivial homology
 groups arise in \v{C}ech filtrations, the space cannot be hyperconvex, but
 only  $\lambda$-hyperconvex for some $\lambda >1$. But  every complete metric
 space is $\lambda$-hyperconvex for some $1\leq\lambda\leq2$,
 c.f. \cite{Grunbaum60}.  (In the discrete case, one might work also with
 $\delta$-hyperbolicity for $\delta>0$.)\\
 From that perspective, hyperconvex spaces are the
 simplest model spaces, and homology can be seen as a topological measure for
 the deviation from such a model. However, this geometric interpretation has
 been dismissed in topological data analysis, by considering the Vietoris-Rips
 filtration instead of \v{C}ech, for the benefit of reducing computational
 complexity. Still, it is possible to infer topological information about a
 space from the Vietoris-Rips filtration, based on Hausmann's
 theorem. However, when one samples a metric space, this depends on how dense
 sample is and the results are accurate only  for small radii. For instance, the Vietoris-Rips complexes of $S^1$ admit holes of dimension larger than $1$ as the radius increases, c.f. \cite{Adamaszek17}. \\

Homology groups, and Betti numbers as integer invariants are fundamental
topological invariants.  Geometry can
provide more refined real valued invariants. And after Riemann \cite{Riemann13,Riemann16}, the
fundamental geometric invariants are curvatures.  In our framework, the
essential geometric content of curvature can be extracted for general metric
spaces. The basic class of model spaces for curvature  is given by the
tripod spaces, a special class containing hyperconvex spaces.   From that
perspective, the geometric content of curvature in the abstract setting
considered here is the deviation from the tripod condition.   Euclidean spaces
only have  a subsidiary role, based on a normalization of curvature that
assigns the value $0$ to them.\\
Considering Euclidean spaces as model spaces is traditionally  justified by
the fact that spaces whose universal cover has  synthetic curvature $\le 0$ in
the sense of Alexandrov are homotopically trivial in the sense that their
higher homotopy groups vanish. In technical terms, they are $K(\pi,1)$ spaces,
with $\pi$ standing for the first homotopy group. The perspective developed
here, however, is a homological and not a homotopical one, and therefore, our
natural comparison spaces are tripods. We have started their investigation in \cite{Joharinad19,Joharinad20}. A more systematic investigation of their properties should be of interest. 
\\
In order to get stronger topological properties, like those of hyperconvex spaces, which are homologically trivial, we might need conditions involving collections of more than three points.\\ 
In fact, according {\cite[Theorem 4.2]{Lindenstrauss62}}, If $X$ is a tripod Banach space on which every collection of four
closed balls $\{B(x_i,r_i)\}_{i=1}^4$ with non-empty pairwise intersection has 
a non-void intersection, then every finite family of closed balls with
non-empty pairwise intersection has also a non-trivial intersection.  In this case, the Vietoris-Rips and \v{C}ech complexes coincide.\\

One can also think about higher order relations and how they can be obtained
from sub-relations (that is from the  relations existing in all subsets of
some smaller size). For instance, in some metric spaces,  a family of $n$
balls has a common point if every subfamily of size $k$ in it has a non-empty
intersection. \cite{Lindenstrauss62} calls this property the $(n,k)$-intersection property. For instance, Helly's theorem
says that Euclidean space $\R^d$ has the $(n,d+1)$-intersection property for
$n\ge d+2$. For a given metric space,  one can compute the deviation from such a property.  \\
From the perspective of \v{C}ech complexes,  this deviation could be quantified by the scaling parameter needed to fill an  $(n-1)$-simplex after all the faces of
dimension $k-1$ are filled. The quantitative measure we introduced provides
us with the scaling function to fill a $2$-simplex after its $1$-dimensional
boundary faces are filled. \\

\bibliographystyle{plain}    
\bibliography{bibtop} 

\end{document}